\documentclass{article}

\usepackage{amssymb}

\usepackage[english]{babel}
\usepackage{amsmath}

\newtheorem{theorem}{Theorem}
\newtheorem{lemma}[theorem]{Lemma}
\newtheorem{e-proposition}[theorem]{Proposition}
\newtheorem{corollary}[theorem]{Corollary}
\newtheorem{e-definition}[theorem]{Definition\rm}

\newcommand{\p}{\mathbb{P}}
\newcommand{\PP}{\mathbb P}

\newcommand{\E}{\mathbb{E}}
\newcommand{\Ex}{\mathbb{E}}

\newcommand{\R}{\mathbb{R}}

\renewcommand{\r}{\right}
\def\la{\left\langle}
\def\ra{\right\rangle}

\def\og{\leavevmode\raise.3ex\hbox{$\scriptscriptstyle\langle\!\langle$~}}
\def\fg{\leavevmode\raise.3ex\hbox{~$\!\scriptscriptstyle\,\rangle\!\rangle$}}

\title{Geometry of log-concave Ensembles of random matrices and approximate reconstruction\thanks{The research was conducted while the authors
participated in the Thematic Program on Asymptotic Geometric Analysis at the Fields Institute in Toronto in Fall 2010.}}

\author{Rados{\l}aw ADAMCZAK\thanks{Research partially supported by MNiSW Grant no. N N201 397437 and the Foundation for Polish Science.}\and
Rafa{\l} LATA{\L}A\thanks{Research partially supported by MNiSW Grant no. N N201 397437 and the Foundation for Polish Science.}\and
Alexander E. LITVAK\and
Alain PAJOR\and
Nicole TOMCZAK-JAEGERMANN\thanks{This author holds the Canada Research Chair in
  Geometric Analysis.}}

\begin{document}
\maketitle

\begin{abstract}
  We study the
  Restricted Isometry Property of a random matrix $\Gamma$ with
  independent isotropic log-concave rows. To this end, we introduce a
  parameter $\Gamma_{k,m}$ that controls uniformly the operator norm
  of sub-matrices with $k$ rows and $m$ columns.  This parameter is
  estimated by means of new tail estimates of order statistics and
  deviation inequalities for norms of projections of an isotropic
  log-concave vector.\\

\end{abstract}

\noindent {\bf Introduction}

The purpose of this short note is to present some new results concerning
geometric properties of random matrices with independent log-concave
isotropic rows obtained recently by the authors. The proofs are deferred
to an upcoming longer article.

Let $T\subset\mathbb{R}^N$ and  $\Gamma$ be an $n\times N$
matrix. Consider the problem of reconstructing any vector $x\in T$
from the data $\Gamma x\in\mathbb{R}^n$, with a fast algorithm.
Clearly one needs some a priori hypothesis on the subset $T$ and of
course, the matrix $\Gamma$ should be suitably chosen. The common
and useful hypothesis is that $T$ consists of sparse vectors, that
is vectors with short support. In that setting, Compressed Sensing
provides a way of reconstructing the original signal $x$ from its
compression $\Gamma x$ with $n\ll N$ by the so-called
$\ell_1$-minimization method. The problem of reconstruction can be
reformulated after D.~Donoho \cite{D1} in a language of high
dimensional geometry, namely, in terms of neighborliness of polytopes obtained by
taking the convex hull of the columns of $\Gamma$.
In this spirit, the sensing matrix is described   by
its columns. From another point of view, the matrix $\Gamma$ may be
also determined by measurements, e.g. by its rows.

Let $0\le m\le N$. Denote by $U_m$ the subset of unit vectors in  $\R^N$,
which  are $m$-sparse, i.e. have at most $m$
non-zero coordinates.  The natural scalar product, the Euclidean
norm and the unit sphere are denoted by $| \la \, \cdot \, ,  \cdot\,  \ra |$,
$|\cdot|$ and $S^{N-1}$.
We also denote by the same notation $|\cdot |$ the cardinality of a set.
For any $x=(x_i)\in\R^n$ we let  $ \|x\|_\infty=\max_{i}|x_i|$. By
$C$, $C_1$, $c$ etc. we will denote absolute positive constants.

Let $\delta_m=\delta_m(\Gamma)=\sup_{x\in U_m} \left| {|\Gamma
x|^2}-\E |\Gamma x|^2 \right |$ be  the Restricted Isometry Property
(RIP) parameter of order $m$. This concept was introduced by
E.~Candes and T.~Tao in \cite{CT1} and  its  important  feature
 is that if $\delta_{2m}$ is appropriately small then
every $m$-sparse vector $x$ can be reconstructed from its
compression $\Gamma x$ by the $\ell_1$-minimization method.
The goal now is to check this property for certain models of matrices.

The articles \cite{ALPT}, \cite{ALPT1} and \cite{cras_alpt}
considered random matrices with independent {\em columns}, and
investigated high dimensional geometric properties of the convex
hull of the columns and the RIP for various models of matrices,
including the log-concave Ensemble build with independent
isotropic log-concave columns. It was shown that various
properties of random vectors can be efficiently studied via operator
norms and the parameter $\Gamma_{n,m}$ recalled below.
In order to control this parameter
an efficient technique of chaining was developed in
\cite{ALPT} and \cite{ALPT1}.

In \cite{MPT}, the authors studied the RIP and more generally the
parameter $\delta _T=\sup_{x\in T} \left| {|\Gamma x|^2}-\E |\Gamma
x|^2 \right |$ for random matrices with independent {\em rows} under
the hypothesis that they are isotropic subgaussian.
It is natural to ask whether random matrices with independent
isotropic log-concave {\it rows} also have the RIP.

Fix integers $ n, N \ge 1$. Let $Y_1, \ldots, Y_n$ be independent
random vectors in  $\R^N$ and let $\Gamma$ be the $n \times N$
random matrix with rows $Y_i$.  Let $T \subset S^{N-1}$ and  $1 \le
k \le n$ and define the parameter  $ \Gamma_k(T)$ by
\begin{equation}
  \label{akT}
  \Gamma_k(T) ^2 = \sup_{y\in T}
\sup_{{I \subset \{1,\ldots,n\}}\atop {|I| = k}}
\sum _{i\in I}  | \la Y_i , y \ra |^2.
\end{equation}

We also denote $\Gamma_{k,m}=\Gamma_k(U_m)$. The role of this
parameter with respect to the RIP is revealed by the following lemma
which reduces a concentration inequality to a deviation inequality.

\begin{lemma}
\label{dva}
Let $Y_1, \ldots, Y_n$ be independent isotropic random vectors in
$\R^N$. Let $T \subset S^{N-1}$ be a finite set.
  Let $0<\theta < 1$ and $B\geq 1$.
Then with probability at least $    1-   |T| \exp\left( - {3 \theta
^2 n}/{8 B^2} \right)$ one has
$$
 \sup_{y\in T} \left|\frac{1}{n} \sum_{i=1}^n(|\langle
 Y_i, y\rangle|^2 - \E |\langle Y_i, y\rangle|^2) \right|
 \leq  \theta + \frac{1}{n}\left( \Gamma_{k}( T) ^2
+ \E  \Gamma_{k}(T)^2\right),
$$
where $k\leq n$ is the largest  integer satisfying
$k\leq (\Gamma_{k}(T)/B)^2$.
\end{lemma}

In this note we focus on the compressed sensing setting where $T$ is
the set of sparse vectors. Lemma~\ref{dva} shows that after a
suitable discretisation, checking the RIP reduces to estimating
$\Gamma_{k,m}$. The idea of such an approach, when $k=n$, originated
from the work of J. Bourgain \cite{B} on the empirical covariance
matrix. It was developed in \cite{ALPT} and \cite{cras_alpt} (with
$T=U_m$), where the estimate of $\Gamma_{n,m}$ played a central role
for solving the Kannan-Lov\'asz-Simonovits conjecture; and   it was
studied  in \cite{M1} where $\Gamma_{k}(T)$ was estimated by means
of Talagrand $\gamma$-functionals.

Using  Lemma \ref{dva} it can be shown  (cf., \cite{cras_alpt} for a
similar argument) that if $0<\theta < 1$, $B\geq 1$, and $m \le
N$ satisfies
$  m \log (C N/ m) \leq 3 \theta ^2 n/16 B^2,$
then with probability at least
$1-  \exp\left( -  {3 \theta ^2 n}/{16 B^2} \right)$ one has
\begin{equation}
  \label{eq:raz}
  \delta_m(\Gamma/\sqrt n)= \sup_{y\in U_m} \left|\frac{1}{n}
\sum_{i=1}^n(|\langle
 Y_i, y\rangle|^2 - \E |\langle Y_i, y\rangle|^2) \right|
 \leq C \theta + \frac{C}{n}\left( \Gamma_{k,m}^2
+ \E \Gamma_{k,m}^2\right),
\end{equation}
where $k\leq n$ is the largest integer satisfying $k\leq
(\Gamma_{k,m}/B)^2$ (note that $k$ is a random variable).

We consider the log-concave Ensemble of $n\times N$ matrices with
independent  isotropic log-concave rows. Recall that a random vector is isotropic log-concave if it is centered,
its covariance matrix is the identity and its distribution has a
log-concave density. Our goal is to bound
 $\Gamma_{k,m}$ for this~Ensemble. This leads  to questions
that require a deeper understanding of some geometric parameters of
log-concave measures, such as tail estimates for order statistics
and deviation inequalities  for  norms of projections.

\smallskip
\noindent{\bf Main results}
\smallskip

Our main theorem provides  upper estimates for
$\Gamma_{k,m}$ valid with large probability for matrices from the
log-concave Ensemble (Theorem \ref{est_akm}). To achieve this we need
some intermediate steps also of a major importance. The first one is a
strengthening of Paouris' theorem (\cite{Pao}) which originally states
that there exists $C >0$ such that
  for every isotropic log-concave vector $X$,
$(\E|X|^p)^{1/p} \le
C((\E|X|^2)^{1/2} + p)$ for $p \ge 1$.
We  define a natural parameter $\sigma_X(p)$ by
$$
\sigma_X(p) = \sup_{t \in S^{N-1}}(\E|\langle t,X\rangle|^p)^{1/p}.
$$
It is known that if $X$ is  isotropic log-concave then
\begin{displaymath}
\sigma_X(p)\leq p\sup_{t\in S^{N-1}}(\Ex|\langle
t,X\rangle|^2)^{1/2}=p.
\end{displaymath}

\begin{theorem}
\label{imprPaouris}
For any $N$-dimensional log-concave vector $X$ and $p\ge 1$ we have
\[
(\Ex|X|^p)^{1/p}\leq C((\Ex |X|^2)^{1/2}+\sigma_X(p)).
\]
\end{theorem}
Our proof of this theorem follows in part the original argument of
\cite{Pao}, which is then complemented by a new analysis of geometry
of log-concave densities.
We were informed by
G.~Paouris that
he has also obtained this result. Another extension is the following
bound on deviations of norm of projections of an isotropic log-concave
vector, uniform over all coordinate projections $P_I$ of a fixed rank.

\begin{theorem}
\label{unifPaouris}
Let $m\leq N$ and $X$ be an isotropic log-concave vector
in $\R^N$. Then for every $t\geq 1$ one has
\[
\p\left(\sup_{{I \subset \{1,\ldots,N\}}\atop {|I| = m}}
|P_IX|\geq Ct\sqrt{m}\log\left(\frac{eN}{m}\r)\r)\leq
\exp\left(-t\frac{\sqrt{m}}{\sqrt{\log (em)}}\log\left(\frac{eN}{m}\r)\r).
\]
\end{theorem}
This theorem is sharp up to $\sqrt{\log (em)}$ in the probability
estimate as the case of a vector with independent exponential
coordinates shows. Actually our further applications require a
stronger result in which the bound for probability is improved by
involving the parameter $\sigma_X$ and its inverse $\sigma_X^{-1}$,
namely

\begin{theorem}
\label{imprunifPaouris}
Let  $m\leq N$ and  $X$ be an isotropic log-concave vector in
$\R^N$. Then for any $t\geq 1$,
\[
\PP\left(\sup_{{I \subset \{1,\ldots,N\}}\atop {|I| = m}}|P_IX|\geq
Ct\sqrt{m}\log\left(\frac{eN}{m}\r)\r)
\leq
\exp\left(-\sigma_{X}^{-1}\left(\frac{t\sqrt{m}\log\left(\frac{eN}{m}\r)}{
\sqrt{\log(em/m_0)}}\r)\r),
\]
where
$
m_0=m_0(X,t)=\sup\left\{k\leq m\colon\ k\log\left(eN/k\r)\leq
\sigma_{X}^{-1}\left(t\sqrt{m}\log\left(eN/m\r)\r)\r\}.
$
\end{theorem}

Theorem \ref{imprunifPaouris} is based on tail estimates for order
statistics of isotropic log-concave vectors.  By $X^*(1)\geq
\ldots\geq X^*(N)$ we denote the non-increasing rearrangement
of $|X(1)|,\ldots,|X(N)|$.  Combining Theorem \ref{imprPaouris} with
methods of \cite{L} we obtain

\begin{theorem}
\label{orderstat} Let $X$ be an $N$-dimensional isotropic log-concave
vector. Then for every $t\geq C\log (eN/\ell)$,
\[
\PP(X^*(\ell)\geq t)\leq \exp(-\sigma_X^{-1}(C^{-1}t\sqrt{\ell})).
\]
\end{theorem}

Introduction of the parameter $\sigma_X$ enables us to obtain
new inequalities for convolutions of log-concave measures. Let
$X_1,\ldots,X_n$ be independent isotropic log-concave random
vectors in $\R^N$. We will consider weighted sums of the vectors
$X_i$ of the form $Y = \sum_{i=1}^n x_iX_i$, where $x =
(x_1,\ldots,x_n) \in \R^n$. Bernstein's inequality and $\psi_1$
estimate for isotropic log-concave random vectors give $\sigma_Y(p)
\le C(\sqrt{p}|x| + p\|x\|_\infty)$ for $p \ge 1$. Together with
Theorem \ref{imprunifPaouris} this yields the following

\begin{corollary}
\label{singlex} Assume that $|x|\leq 1$ and  $1\geq b\geq
\max(\|x\|_{\infty},1/\sqrt{m})$. Then for any $t\geq 1$,
\[
\PP\Bigg(\sup_{{I \subset \{1,\ldots,N\}}\atop {|I| = m}}|P_IY|\geq
Ct\sqrt{m}\log\Big(\frac{eN}{m}\Big)\Bigg)
\leq
\exp\Bigg(-\frac{t\sqrt{m}
\log\Big(\frac{eN}{m}\Big)}{b\sqrt{\log(e^2b^2m)}}\Bigg).
\]
\end{corollary}

We now pass to bounds on deviation of $\Gamma_{k,m}$.
To get a slightly simplified formula we assume that~$N \ge n$.

\begin{theorem}
\label{est_akm} Let $1 \le n \le N$, and let $\Gamma$ be an $n\times
N$ random matrix with independent isotropic log-concave rows. For
any integers $k\le n$, $m\le N$ and any $t \ge 1$, we have
$$
\p(\Gamma_{k,m} \ge Ct \lambda ) \le \exp(-t
\lambda/\sqrt{\log (3m)}),
$$
where $\lambda = \sqrt{\log\log (3m)} \sqrt{m}\log(e N /m)
   + \sqrt{k}\log(en/k) $.
\end{theorem}

\medskip

The threshold value $\lambda$ in the above theorem is optimal,
up to the factor of $\sqrt{\log\log
  (3m)}$. Assuming additionally
unconditionality of the distributions of the rows, we can remove
this factor and get a sharp estimate.

The proof of the above theorem is composed of two parts, depending on
the relation between $k$ and the quantity
$
  k' = \inf\{\ell \ge 1\colon m\log (eN/m )\le
  \ell\log (en/\ell )\}.
$
First, we adjust the chaining argument from \cite{ALPT} to reduce the
problem to the case $k \le k'$. This step also involves
Theorem~\ref{unifPaouris}. Next, we use Corollary~\ref{singlex}
combined with another chaining to complete the argument.

Theorem \ref{est_akm} together with (\ref{eq:raz}) allows us to
prove the RIP result for matrices $\Gamma$ with independent
isotropic log-concave rows. The result is optimal, up to the factor
$\log\log 3m$, as shown in  \cite{ALPT1}. As for Theorem
\ref{est_akm}, assuming unconditionality of the distributions of the
rows, we can remove this factor.

\begin{theorem}
\label{rip} Let $0<\theta < 1$, $1 \le n \le N$. Let $\Gamma$ be an
$n\times N$ random matrix with independent isotropic log-concave
rows. There exists $c(\theta)>0$ such that
$\delta _m (\Gamma/\sqrt n) \leq \theta$ with overwhelming probability
 whenever
$$m\log^2(2N/n) \log\log 3m \leq c(\theta) n.$$
\end{theorem}

\begin{flushleft}
Rados{\l}aw Adamczak\\
Institute of Mathematics,\\
University of Warsaw\\
Banacha 2, 02-097 Warszawa, Poland \\
{\tt R.Adamczak@mimuw.edu.pl}

\smallskip
Rafa{\l} Lata{\l}a\\
Institute of Mathematics,\\
University of Warsaw\\
Banacha 2, 02-097 Warszawa, Poland \\
{\tt rlatala@mimuw.edu.pl}

\smallskip
Alexander E. Litvak\\
Department of Mathematical and Statistical Sciences,\\
University of Alberta,\\
Edmonton, Alberta, Canada T6G 2G1\\
{\tt alexandr@math.ualberta.ca}

\smallskip
Alain Pajor\\
Equipe d'Analyse et Math\'ematiques Appliqu\'ees,\\
Universit\'e Paris Est, \\
5 boulevard Descartes, Champs sur Marne, 77454 Marne-la-Vallee, \\
Cedex 2, France\\
{\tt alain.pajor@univ-mlv.fr}

\smallskip
Nicole Tomczak-Jaegermann,\\
Department of Mathematical and Statistical Sciences,\\
University of Alberta,\\
Edmonton, Alberta, Canada T6G 2G1\\
{\tt nicole@ellpspace.math.ualberta.ca}
\end{flushleft}

\end{document}